\documentclass[a4paper,11pt]{amsart}
\usepackage{amsfonts}
\usepackage{amsthm}
\usepackage{amssymb}
\usepackage{amsmath}
\usepackage{comment}
\usepackage{a4wide}
\usepackage{dsfont}
\usepackage{bbm}
\usepackage{xparse,mathrsfs,mathtools}
\usepackage{enumitem}
\usepackage{color}
\usepackage{caption}
\usepackage{subcaption}
\usepackage{hyperref}
\usepackage{fancyhdr}

\newtheorem{theorem}{Theorem}[section]
\newtheorem{lemma}{Lemma}[section]
\newtheorem{problem}{Problem}[section]
\newtheorem{conjecture}{Conjecture}[section]
\newtheorem{question}{Question}[section]
\newtheorem{remark}{Remark}[section]

\allowdisplaybreaks

\begin{document}

\title[]{Open Problems UP24}
\author[]{Maryna Manskova}

\begin{abstract} 
The conference \textit{Unexpected Phenomena in Energy Minimization and Polarization}, held in Sofia, Bulgaria in 2024, provided a platform for researchers to discuss and propose challenging open questions across various fields, such as potential theory, approximation, special functions, point configurations, lattices, and numerical analysis. The open problems sessions were productive, fruitful and led to a range of interesting questions. In this document, we present these open problems. 
\end{abstract}
	
\maketitle

\section{\textbf{When is the equilibrium support a sphere?}\\Robert Womersley and Edward Saff}

\, \\
The \textit{Riesz} $s$\textit{-kernel} $K_s\colon \mathbb{R}^d\rightarrow (-\infty,+\infty] $ is defined by
$$K_s(x)=\begin{cases}
    \frac{1}{s\|x\|^s} & \text{if} \ s\neq 0
    \\
    -\log(\|x\|) & \text{if}\ s=0
\end{cases},$$
where $\|x\|=\sqrt{x_1^2+\cdots+x_d^2}$ is the Euclidean norm. We assume $-2<s<d$, which ensures that the kernel is integrable and
conditionally strictly positive definite on compact sets.

Let $V\colon \mathbb{R}^d\rightarrow (-\infty,+\infty] $ be an external field and 
$\mathcal{P}(\mathbb{R}^d)$ be the set of probability measures on $\mathbb{R}^d$. For $s<d$, the \textit{energy} of $\mu\in\mathcal{P}(\mathbb{R}^d)$ is defined by
$$I_{s,V}(\mu)=\iint\left(K_s(x-y)+V(x)+V(y)\right)d\mu(x)d\mu(y)\in (-\infty,+\infty].$$
When they exist, the minimizers, called \textit{equilibrium} measures, are denoted by $\mu_{\text{eq}}$. 

We are particularly interested in radial external fields, i.e. fields of the form
$$V(x)=v(r^2), \ \ r=\|x\|,$$
where $v\colon [0, +\infty) \rightarrow (-\infty, +\infty]$ is lower semi-continuous, bounded from below, and finite
on some interval $(a,b)$. For example,
$$v(\rho)=\gamma \rho^{\alpha/2}\ \text{and}\ V(x)=\gamma \|x\|^\alpha.$$
It was shown in \cite{Chafa2024} that under some conditions on $s$, $d$, and $\alpha$ the support of $\mu_{\text{eq}}$ is a sphere $\mathbb{S}_R^{d-1}$, when radial symmetry implies that the equilibrium measure is the uniform measure $\sigma_R$ on $\mathbb{S}_R^{d-1}$. Let $c_{s,d}$ and $b_d$ be defined by equations (1.4) and (1.5) of \cite{Chafa2024}.
\begin{theorem}\label{Theorem_sphere}
    Suppose that $-2<s<d-3$ and $V(x)=\frac{\gamma}{\alpha} \|x\|^\alpha$, where $\gamma>0$ and $\alpha>\max\{-s,0\}$. Define 
    $$\alpha_{s,d}=\begin{cases}
        \max\left\{\frac{sc_{s,d}}{2-2c_{s,d}},2-\frac{(s+2)(d-s-4)}{2(d-s-3)}\right\} & s\neq 0
        \\
        \max\{-\frac{1}{2b_d},2-\frac{(d-4)}{(d-3)}\} & s=0
    \end{cases}.$$
If $\alpha\geq \alpha_{s,d}$, then $\mu_{\text{eq}}=\sigma_{R_*}$, where
$$R_*=\left(\frac{c_{s,d}}{2\gamma}\right)^{\frac{1}{\alpha+s}}.$$
    Furthermore, the threshold $\alpha_{s,d}$ is a sharp bound for $\alpha$, meaning that if $\max\{-s,0\} < \alpha <
\alpha_{s,d}$, then for all $R>0$, $\sigma_R$ is not a minimizer of $I_{s,V}$. 
\end{theorem}

\begin{figure}[h]
\centering
\includegraphics[width=0.85\textwidth]{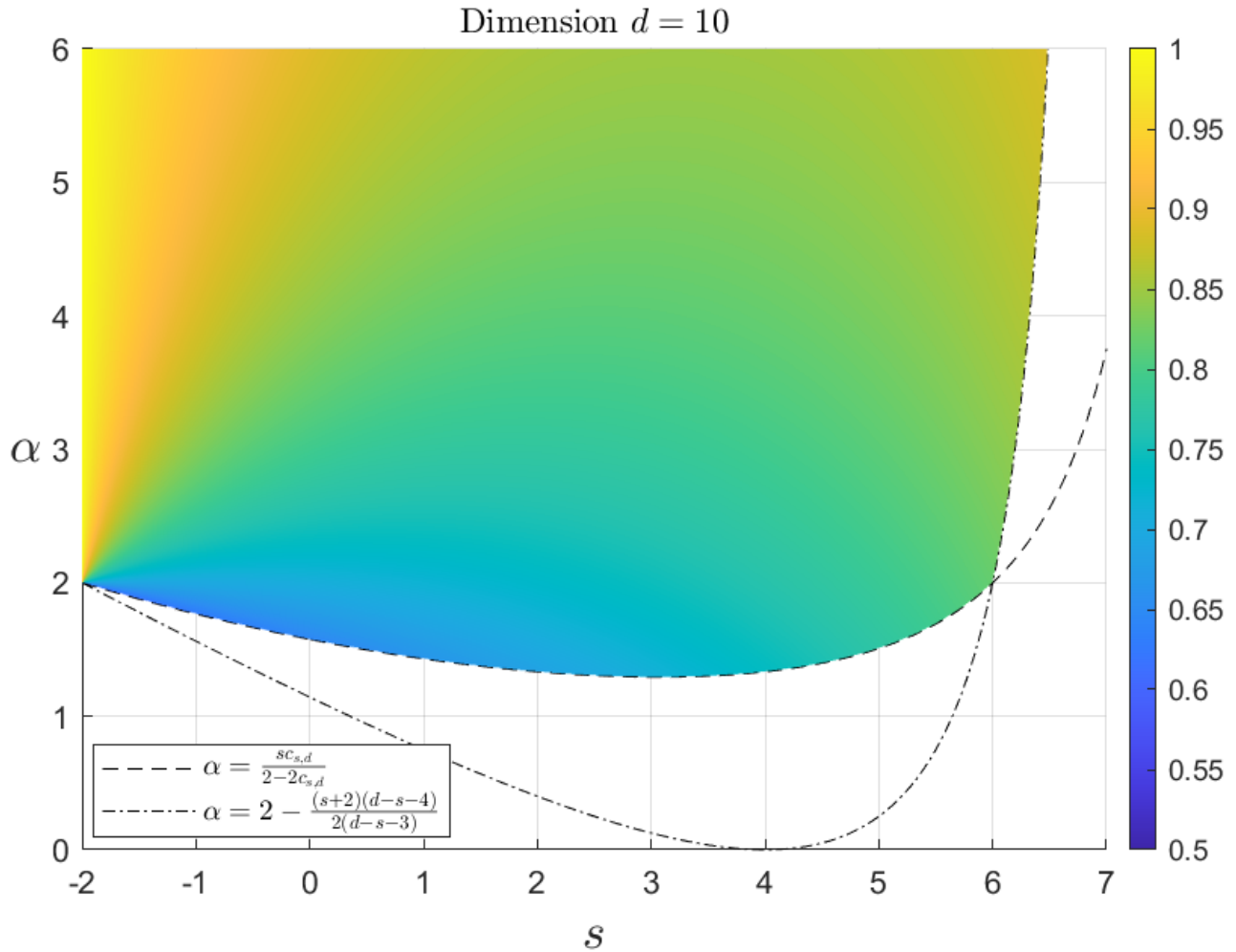}
\caption{ Plot for $d=10$ and Riesz parameter $-2<s<d-3$, where the colour gives the value of $R_*$ when
the equilibrium support is $\mathbb{S}^{d-1}_{R_*}$. The external field power $\alpha\geq \alpha_{s,d}$ as in Theorem \ref{Theorem_sphere} with $\gamma=1$. Outside the coloured region the support is not a sphere.}
\label{threshold}
\end{figure}

Theorem \ref{Theorem_sphere} provides a characterization of when the support of the
equilibrium measure is a sphere, which still leaves open the following questions.
\begin{question}[Robert Womersley, Edward Saff]
    What can we say about the equilibrium measures for values of $s$ with $d-3<s<d$?
\end{question}
The support of $\mu_{\text{eq}}$ may be a ball, not a sphere. This was shown for the cases $s=d-1$ by M. Riesz (\cite{riesz1931certaines}, \cite{Riesz1988IntegralesDR}) and $s=d-3$ by D. Chafaï, E. B. Saff, R. S. Womersley \cite{Chafa2022}.

In all previous examples, the dimension of  $\text{supp}(\mu_{\text{eq}})$ is less than $d$. The question is whether dimension reduction occurs for other combinations of $-2<s<d$ and $\alpha>\max\{0,-s\}$.
\begin{question}[Robert Womersley, Edward Saff]
    When is the support of the equilibrium measure full dimensional and when is the support of dimension less than $d$? 
\end{question}
 In the case $s=d-4$ and $\alpha<2$, $\mu_{\text{eq}}$  has a full dimensional component (see \cite{Chafa2023}).

The result of Theorem \ref{Theorem_sphere} provides a sharp bound for $\alpha$. It would be interesting to study the case $\alpha<\alpha_{s,d}$ (white region on Fig. \ref{threshold}).
\begin{question}[Robert Womersley, Edward Saff]
    What can we say about the equilibrium measure for values of $\alpha$ with $\max\{-s,0\}<\alpha<\alpha_{s,d}?$
\end{question}

\begin{figure}[h]
\centering
\includegraphics[width=0.80\textwidth]{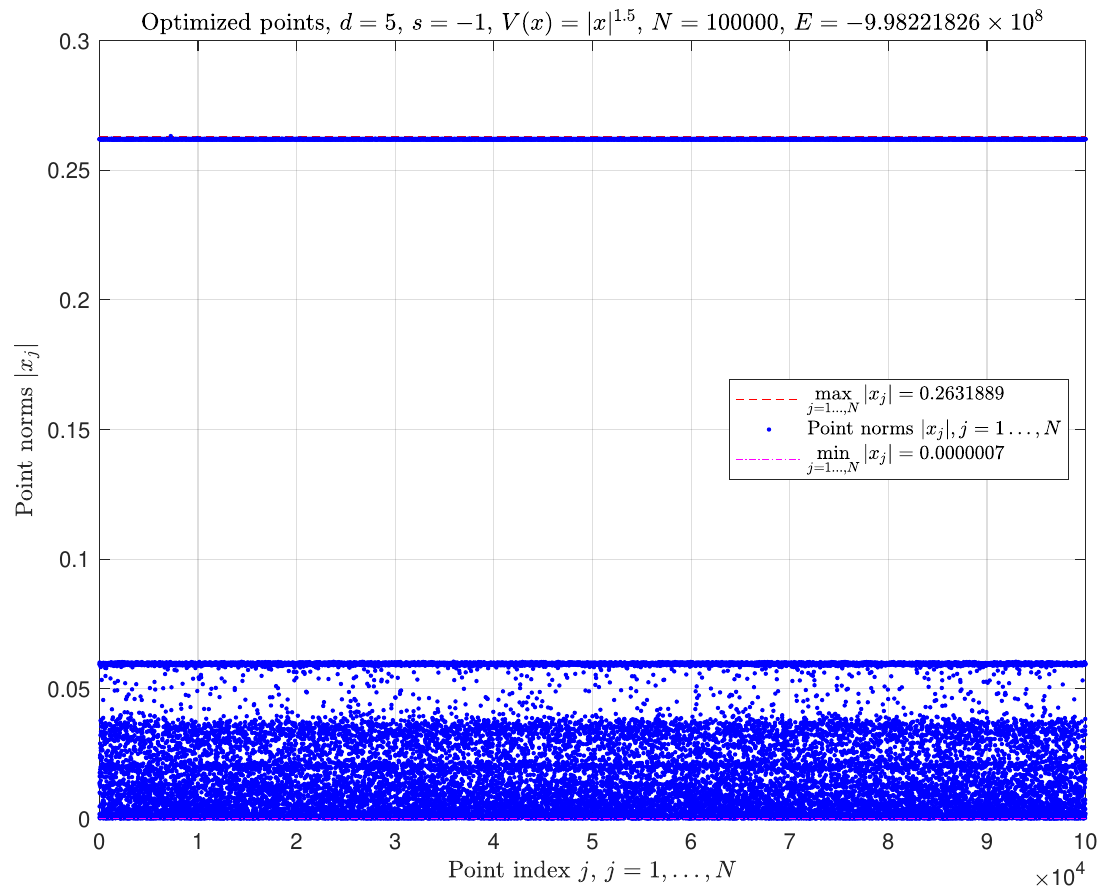}
\caption{ Numerical results for $d=5$, $s=d-6=-1$, and $\alpha=\frac{3}{2}<\alpha_{s,d}$.}
  \label{numerical_result}
\end{figure}

One can get unexpected answers to this question. For example, optimization of the energy for the case $d=5$, $s=-1$, and $\alpha=\frac{3}{2}$ gives us the results in Fig. \ref{numerical_result}. It looks like an outer sphere, an inner sphere and points around the origin. 

\clearpage

\section{\textbf{Inner product expansions for logarithmic and Riesz $\mathbf{s}$-energy}\\Johann Brauchart}

Let us consider the \textit{discrete minimal logarithmic} and \textit{Riesz energy} problem on the unit sphere $\mathbb{S}^d\subset\mathbb{R}^{d+1}$ (w.r.t. to the Euclidean distance). Namely, we are interested in the asymptotic expansion of
$$\mathcal{E}_{\log}(\mathbb{S}^d;N)=\min\limits_{\mathbf{x}_1,\dots,\mathbf{x}_N\in\mathbb{S}^d}\sum\limits_{i\neq j}\log\frac{1}{\|\mathbf{x}_i-\mathbf{x}_j\|}$$
and for $s>0$
$$\mathcal{E}_{s}(\mathbb{S}^d;N)=\min\limits_{\mathbf{x}_1,\dots,\mathbf{x}_N\in\mathbb{S}^d}\sum\limits_{i\neq j}\frac{1}{\|\mathbf{x}_i-\mathbf{x}_j\|^s}\quad \ \  \text{as}\  N\rightarrow\infty.$$

For $d=2$, it is shown in \cite{Btermin2016} that there exists $C_{\log}\neq 0$ independent of $N$ such that
$$\mathcal{E}_{\log}(\mathbb{S}^2;N)=\left(\frac{1}{2}-\log 2\right)N^2-\frac{1}{2}N\log N+C_{\log}N+o(N)\quad\quad \text{as }\  N\rightarrow\infty.$$
Moreover, it turns out that $C_{\text{BHS}}$, earlier conjectured in \cite{Brauchart2012}, is an upper bound for $C_{\log}$.
\begin{align*}
    C_{\log}\leq C_{\text{BHS}}&=2\log 2 +\frac{1}{2}\log\frac{2}{3}+3\log \frac{\sqrt{\pi}}{\Gamma(\frac{1}{3})}
    \\ &=-0.05560530494339251850\dots
\end{align*}

The best known lower bound is given in \cite{Lauritsen2021}, 
\begin{align}\label{log_constant_bound}
    C_{\log}=2e_{\text{Jel}}+\frac{\log 4\pi}{2}\geq \log 2-\frac{3}{4}=-0.05685281944005469\dots,
\end{align}
where $e_{\text{Jel}}$ is the thermodynamic limit of the properly scaled jellium energy of $N$ particles in a domain of measure $N$. An alternative proof of (\ref{log_constant_bound}) and a generalization of this proof for the Green energy is given in \cite{Beltrn2023}. 

One possible approach to finding a lower bound is to expand $\log\frac{1}{\|\mathbf{x}-\mathbf{y}\|}$ as a power series of the inner product $\langle \mathbf{x},\mathbf{y}\rangle$ (see \cite{https://doi.org/10.48550/arxiv.2502.15984}). Using this method, we can recover Steinerberger's lower bound (see \cite{Steinerberger2022}).
\begin{lemma}[Johann Brauchart]
    Let $d=2$. Then 
    $$\mathcal{E}_{\log}(\mathbb{S}^2;N)\geq \left(\frac{1}{2}-\log 2\right)N^2-\frac{1}{2}N\log N+F(c)N+\mathcal{O}(1), $$
    where the function 
    $$F(c)=\frac{\log 2}{2}-\frac{\gamma}{2}-\frac{\log c}{2}-\frac{1}{4c}$$
    has a unique maximum at $c=c^*=\frac{1}{2}$ in $(0,\infty)$ with value 
    $$F(c^*)=\log 2-\frac{1}{2}-\frac{\gamma}{2}=-0.09546065189082112\dots .$$
\end{lemma}
In order to improve the result, we need to find non-trivial bounds for some non-negative summands that were ignored. By applying summation by parts, we can obtain a connection between the minimal Riesz $2$-energy and the logarithmic energy,
\clearpage
\begin{align*}
    \mathcal{E}_{\log}(\mathbb{S}^2;N)\geq &\left(\frac{1}{2}-\log 2\right)N^2-\frac{1}{2}N\log N
    \\&+\left(\log 2-\gamma+\frac{2}{N^2}\left( \mathcal{E}_{2}(\mathbb{S}^2;N)-\frac{1}{4}N^2\log N\right)\right)N
    \\ & +\frac{1}{2}\sum\limits_{m=1}^{M-1}\frac{1}{\ell(\ell+1)}\sum\limits_{m=1}^\ell S_m(\mathbf{x}_1^*,\dots,\mathbf{x}_N^*)-\sum\limits_{m=M+1}^\infty \frac{1}{m(m-1)}\sum\limits_{j\neq k}\frac{\langle \mathbf{x}_j^*,\mathbf{x}_k^*\rangle^m}{2(1-\langle \mathbf{x}_j^*,\mathbf{x}_k^*\rangle)}\\
    &+\mathcal{O}(1)+\mathcal{O}\left(\frac{2}{N^2}\left( \mathcal{E}_{2}(\mathbb{S}^2;N)-\frac{1}{4}N^2\log N\right)\right).
\end{align*}
We recall the conjecture for the expansion of the Riesz $2$-energy.
\begin{conjecture}
    $$\mathcal{E}_2(\mathbb{S}^2;N)=\frac{1}{4}N^2\log N+C_{2,2}N^2+o(N^2),$$
    where \begin{align}\label{C_2,2}
        C_{2,2}=\frac{1}{4}\left(\gamma-\log(2\sqrt{3}\pi)\right)+\frac{\sqrt{3}}{4\pi}\left(\gamma_1(\tfrac{2}{3})-\gamma_1(\tfrac{1}{3})\right)
    \end{align}
and $\gamma_n(\alpha)$ is the generalized Stieltjes constant appearing as the coefficient $\frac{\gamma_n(\alpha)}{n!}$ of $(1-s)^n$ on the Laurent series expansion of the Hurwitz zeta function $\zeta(s,a)$ about $s=1$.
\end{conjecture}
The connection formula also provides the following relation between coefficients of the expansions of the energies.
\begin{theorem}[Johann Brauchart]
    $$C_{\text{BHS}}-2C_{2,2}=\log 2-\gamma>0,$$
where $C_{2,2}$ is defined by (\ref{C_2,2}).
\end{theorem}
\begin{question}[Johann Brauchart]
    Can other relation be established between coefficients of different (w.r.t energy) asymptotic expansions?
\end{question}

\begin{question}[Johann Brauchart]
Are there other more efficient resummation methods?
\end{question}

\begin{question}[Johann Brauchart]
Are there other replacements for inner product expansion?
\end{question} 

\begin{remark}
The ideas and arguments presented in \cite{https://doi.org/10.48550/arxiv.2502.15984} for computing the spherical cap $\mathbb{L}_2$-discrepancy (a continuous kernel) have been extended to logarithmic and Riesz energies (singular kernels) in \cite{Paper2}.
\end{remark}

\clearpage

\section{\textbf{Energy minimization on convex bodies}\\ Ryan Matzke}
Let $C\subset \mathbb{R}^d$ be a compact convex body. For $s<d$, the \textit{energy} of a probability measure  $\mu$ on $C$ is defined by
$$I_{s}(\mu)=\int\limits_C\int\limits_C K_s(x,y)d\mu(x)d\mu(y),$$

where $K_s(x,y)$ is the Riesz $s$-kernel. The minimizers of $I_s$ are called \textit{equilibrium} measures and denoted by $\mu_{\text{eq}}$. 

In general case, it is known that the minimizers of the energy integral can be characterized by:

1) $d-2<s<d$, \quad   $\text{supp}(\mu_{eq})=C$;

2) $s\leq d-2$, \quad $\text{supp}(\mu_{eq})\subseteq \partial C$;

3) $s<-1$, \quad $\text{supp}(\mu_{eq})$ is a subset of extreme points;

4) $s<-2$, \quad $\text{supp}(\mu_{eq})$ is discrete with at most $d+1$ points.
\\
\\
In the case of a ball $B_d(0,R)$, we know that the equilibrium measures are given by

1) $d-2<s<d$, \quad   $d\mu_{eq}(x)=c(R^2-\|x\|^2)^{-\frac{d-s}{2}}dx$;

2) $-2<s\leq d-2$, \quad $\mu_{eq}=\sigma$ is the uniform measure on $\partial B$;

3) $s<-2$, \quad $\mu_{eq}=\frac{1}{2}(\delta_p+\delta_{-p})$, where $p\in\partial B$.
\\
\\
At the same time, for a hypercube $Q_d$ in $\mathbb{R}^d$, the support of $\mu_{eq}$ consists of the vertices of $Q_d$ for $-2<s\leq-1$. It looks like the fact that $B_d(0,R)$ and $Q_d$ have different structure (namely, $Q_b$ has lower dimensional faces) enables us to hit lower dimensional (discrete) minimizers earlier. 

For a larger interval $-2<s\leq q-2$ for some $q\in\{1,\dots,d\}$ it is known that $\text{supp}(\mu_{\text{eq}})$ is contained in the union of $(q-1)$- dimensional faces of $Q_d$. This gives us the following bounds for the dimension of $\mu_{\text{eq}}$,
\begin{align}\label{dim}
    s\leq \dim(\text{supp}(\mu_{\text{eq}}))\leq q-1.
\end{align}
\\
Let $s=0$ and $q=2$, then $0\leq \dim(\text{supp}(\mu_{\text{eq}}))\leq 1.$ 
\\
\begin{problem}[Ryan Matzke]
    Improve lower bound in (\ref{dim}).
\end{problem}

Even if the dimension is known, there is still a question about the structure of the support of minimizers.
\begin{problem}[Ryan Matzke]
  Find more information about the support of $\mu_{eq}$.
\end{problem}

Written above can be formulated not only for hypercubes but also for any convex polytope.  The question is if one can get similar results for a larger class of convex bodies. A good example could be the Reuleaux triangle.

\begin{question}[Ryan Matzke]
    Can we obtain similar results for other convex bodies with lower dimensional faces?
\end{question}

\clearpage

    \section{\textbf{Minimizers of the $\mathbf{p}$-frame energy. Distances on the torus.}
    \\ Dmitriy Bilyk}

Let $F(t)=|t|^p$ with $p>0$. We study minimizers of the \textit{$p$-frame energy} 
$$I_F(\mu)=\int\limits_{\mathbb{S}^{d-1}}\int\limits_{\mathbb{S}^{d-1}}|\langle x,y\rangle|^p d\mu(x)d\mu(y).$$

In the case $p=2$, we have tight frames and isotropic measures on $\mathbb{S}^{d-1}$ (see \cite{Benedetto2003}). For $0<p<2$ the minimizers are orthonormal bases (see \cite{Ehler2012}). For all even $p>2$ we obtain spherical designs and the uniform probability measure $\sigma$. In general, for all other values of $p$ the problem is not solved. Numerical experiments showed that the minimizers are discrete. Moreover, it was proved that if $p$ is not a even integer, then supports of minimizers have empty interior. The next step would be to find some bounds on the dimension of supports. 
\begin{conjecture}[D. Bilyk, A. Glazyrin, R. Matzke, J. Park, O. Vlasiuk, \cite{Bilyk2021,Bilyk2022}]
    For $p>2$ and $p\neq 2k$, $k\in\mathbb{N}$, all minimizers of the $p$-frame energy are discrete measures. 
\end{conjecture}

For $p$ such that  $\sigma$ is not a minimizer, there are negative coefficients in the Gegenbauer series expansions of $|t|^p$. As $p$ grows, these negative coefficients start later and later. 

\begin{problem}[Dmitriy Bilyk]
    Let $\kappa\colon [-1,1]\rightarrow\mathbb{R}$. Assume that there exist negative coefficients $c_l$ in the Gegenbauer expansion of  $\kappa(t)$, where $l>N$ for some large $N$. If there exists a discrete minimizer $\mu^*$ of $I_{\kappa}$, what can one say about the size of the support of $\mu^*$?
\end{problem}

Let $\alpha>0$ and $\Omega\subset\mathbb{R}^n$ be compact. For a given distance function $\rho$, we maximize the energy integral
$$I_\alpha(\mu)=\int\limits_{\Omega}\int\limits_{\Omega}\rho(x,y)^\alpha d\mu(x)d\mu(y).$$

\begin{figure}[h]
     \centering
     \begin{subfigure}[b]{0.49\textwidth}
         \centering
         \includegraphics[width=\textwidth]{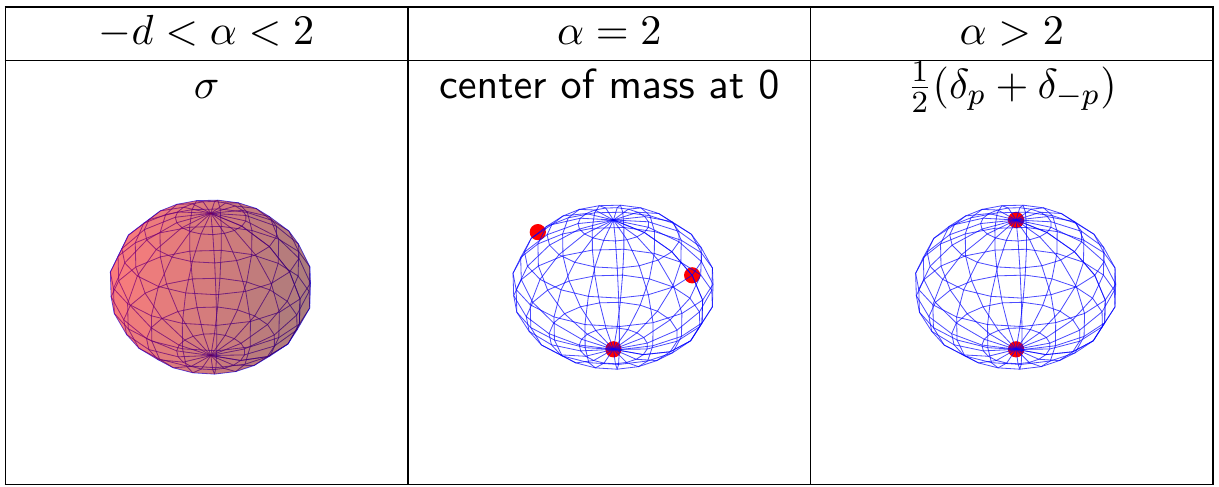}
         \caption{Euclidean distance}
     \end{subfigure}
\hfill
     \begin{subfigure}[b]{0.49\textwidth}
         \centering
         \includegraphics[width=\textwidth]{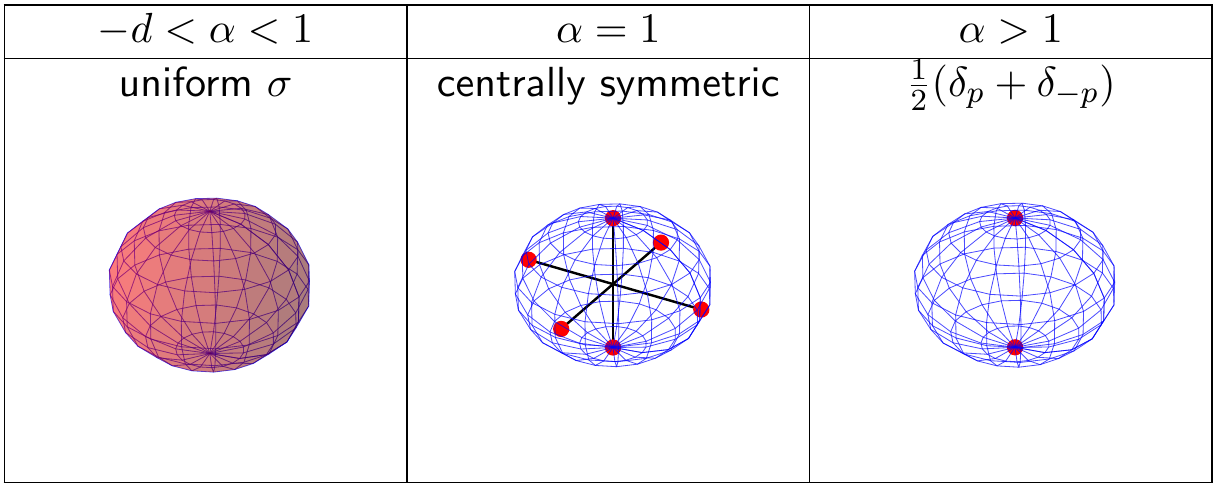}
         \caption{Geodesic distance}
     \end{subfigure}
     
        \caption{Maximizers of $I_\alpha$.}
        \label{fig:spheres}
\end{figure}

In the case $\Omega=\mathbb{S}^d\subset \mathbb{R}^{d+1}$, we have results for the Euclidean and the geodesic distances (see \cite{Bjrck1956} and \cite{Bilyk2018}). The phase transition points between $\sigma$ and $\frac{1}{2}(\delta_p+\delta_{-p})$ are $\alpha=2$ and $\alpha=1$ respectively (see Fig. \ref{fig:spheres}).

Another example that has features of both Euclidean and geodesic distances is the flat torus. We know that $\sigma$ maximizes $I_\alpha$ for $\alpha=1$. In the case $\alpha=2$, maximizers are $\frac{1}{2}(\delta_p+\delta_{p^*})$, where $p^*=p+(\frac{1}{2},\frac{1}{2})$. It means that the phase transition point(s) must be somewhere in the interval $(1,2)$.
\begin{problem}[Dmitriy Bilyk]
    Let $\Omega=\mathbb{T}^2$ be the flat torus with $\rho(x,y)=\min\limits_{k\in\mathbb{Z}^2}\|x-y+k\|$. What are the phase transition points in this case?
\end{problem}

\clearpage

    \section{\textbf{Riesz energy for geodesic and chordal distances}\\ Peter Grabner}
Let $\Omega$ be a sphere $\mathbb{S}^d$ or a projective space $\mathbb{FP}^d$ equipped with a distance $\rho(x,y)$ equivalent to the geodesic distance. For $-2<s<d$, we consider the \textit{Riesz energy} 
$$I_{s}(\mu)=\int\limits_\Omega\int\limits_\Omega K_s(\rho(x,y))d\mu(x)d\mu(y).$$
Let us denote by $S(\Omega)$ the set of $s$ such that the uniform probability measure $\sigma$ on $\Omega$ is a minimizer for $I_s$. We want to study the structure of the set $S(\Omega)$. It is known that if $\sigma$ is a minimizer of $I_{s^*}$ for some $s^*\leq 0 $, then $(s^*,d)\subset S(\Omega)$ (see \cite{Bilyk_Grabner2025:positive_definite_singular}). In all known examples, the set $S$ is connected. Hence it makes sense to ask if the statement holds for $s^*>0$.
\begin{question}[Peter Grabner]
   Is the following statement true: let $s^*>0$ such that $\sigma$ is minimizer of $I_{s^*}$, then it is a minimizer for all $s^*<s<d$.
\end{question}

\clearpage

\section{\textbf{Integration weights}\\Jordi Marzo}
Let us denote by $\mathbb{P}_n^r$ the space of spherical polynomials on $\mathbb{S}^r$ of degree at most $n$, and $N=\dim \mathbb{P}^r_n$. Consider a fundamental system of points $\{x_1,x_2,\dots,x_N\}\subset \mathbb{S}^r$, i.e.  the determinant of the interpolation matrix is non-zero. We can build the fundamental Lagrange polynomials $\{\ell_i(x)\}_{i=1}^N$ associated with the given fundamental system, this can be written as
$$ \ell_i(x)\in \mathbb{P}_n^r, \quad  \ell_i(x_j)=\delta_{ij},\ \ \  i,j = 1,\dots,n. $$
The \textit{integration weights} are defined by 
$$w_i=\int_{\mathbb{S}^r_n}\ell_i(x)d\sigma(x).$$

In the case $n=2$, it was shown in \cite{Reimer1994} that sets of Fekete points give positive integration weights; moreover, they are close to $1/N$. For larger values of $n$ (from $1$ to $51$, then for $56$, $63$, $64$, $72$, $96$, $127$, $128$, and $191$), the positivity of the weights was shown numerically in \cite{Sloan2004}.

One can try to determine the behaviour of the (random) integration weights of some other (random) configurations of points.

\begin{question}[Jordi Marzo]
    
Consider random points on the sphere (uniform i.i.d or from a determinantal point process). Can one quantify how good the points are using the notion of positivity of the integration weights (or being close to $1/N$)?
\end{question}

\clearpage

\section{\textbf{Snake Polynomials}\\Geno Nikolov}

    Let us denote by $\mathcal{P}_m$ be the set of all algebraic polynomials of degree at most $m$. Without loss of generality we can assume that the polynomials have real coefficients.

Take a \textit{majorant} $\mu\in C[-1,1]$, $\mu(x)\geq 0$, $x\in[-1,1]$. If there exists a non-zero polynomial $P\in \mathcal{P}_n$ such that $-\mu(x)\leq P(x)\leq \mu(x),$ $x\in[-1,1]$, then there exists a unique (up to orientation) polynomial $\omega_\mu\in\mathcal{P}_n$ which oscillates most between $\pm\mu$. We call $\omega_\mu$ a \textit{snake polynomial} associated with $\mu$.

The $n$-th snake polynomial $\omega_\mu$ is uniquely determined by the following properties:

a) $|\omega_\mu(x)|\leq \mu(x)$ for all $x\in[-1,1]$;

b) There exists a set $\delta^*=(\tau^*_i)^n_{i=0}$, $-1\leq \tau_n^*<\cdots <\tau^*_0\leq 1$ such that 
$$\omega_\mu(\tau_i^*)=(-1)^i\mu(\tau_i^*),\quad i=0,\dots,n,$$
where $\delta^*$ is referred to as the set of \textit{alternation points} of $\omega_\mu$.
\begin{problem}[Geno Nikolov]
Find a class of majorants $\mu$ for which the associated snake polynomials have non-negative (or sign-alternating) expansion in the Chebyshev polynomials of the first kind. 
\end{problem}

\begin{figure}[h]
     \centering
     \begin{subfigure}[b]{0.3\textwidth}
         \centering
         \includegraphics[width=\textwidth]{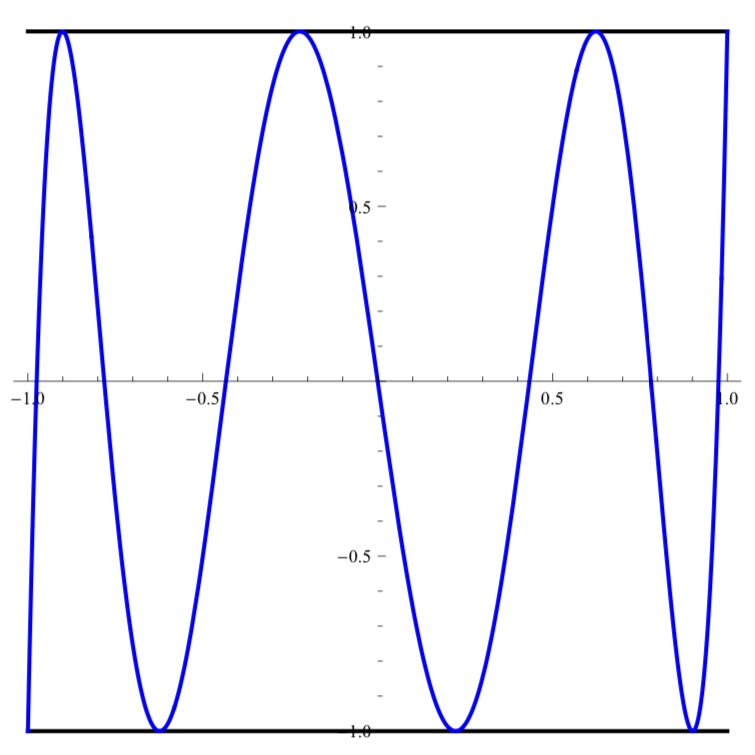}
         \caption{$\mu(x)=1$}
         \label{fig:1}
     \end{subfigure}
     \hfill
     \begin{subfigure}[b]{0.3\textwidth}
         \centering
         \includegraphics[width=\textwidth]{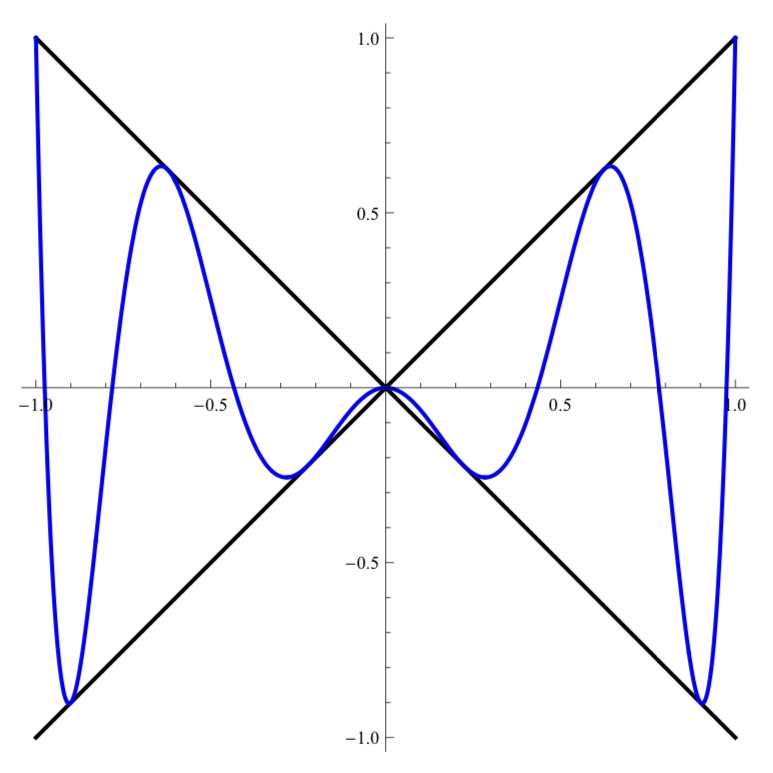}
         \caption{$\mu(x)=|x|$}
         \label{fig:2}
     \end{subfigure}
     \hfill
     \begin{subfigure}[b]{0.3\textwidth}
         \centering
         \includegraphics[width=\textwidth]{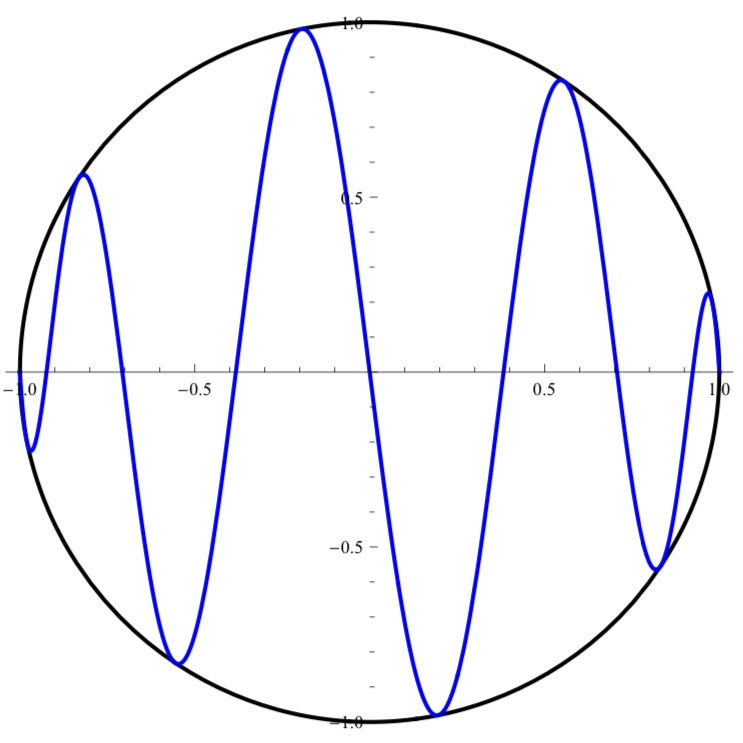}
         \caption{$\mu(x)=\sqrt{1-x^2}$}
         \label{fig:3}
     \end{subfigure}
     \\
     \begin{subfigure}[b]{0.3\textwidth}
         \centering
         \includegraphics[width=\textwidth]{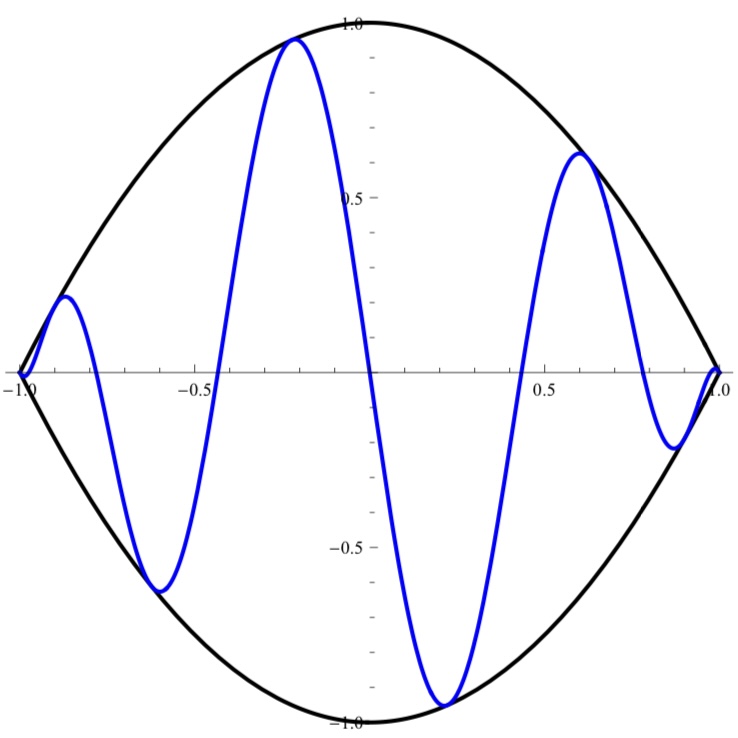}
         \caption{$\mu(x)=1-x^2$}
         \label{fig:4}
     \end{subfigure}
     \hfill
     \begin{subfigure}[b]{0.3\textwidth}
         \centering
         \includegraphics[width=\textwidth]{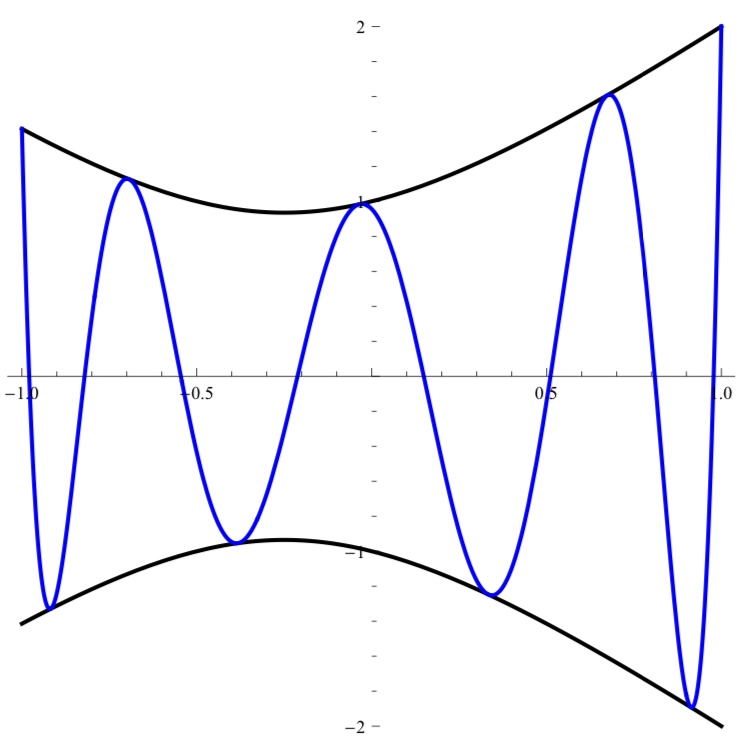}
         \caption{$\mu(x)=\sqrt{2x^2+x+1}$}
         \label{fig:5}
     \end{subfigure}
     \hfill
     \begin{subfigure}[b]{0.3\textwidth}
         \centering
         \includegraphics[width=\textwidth]{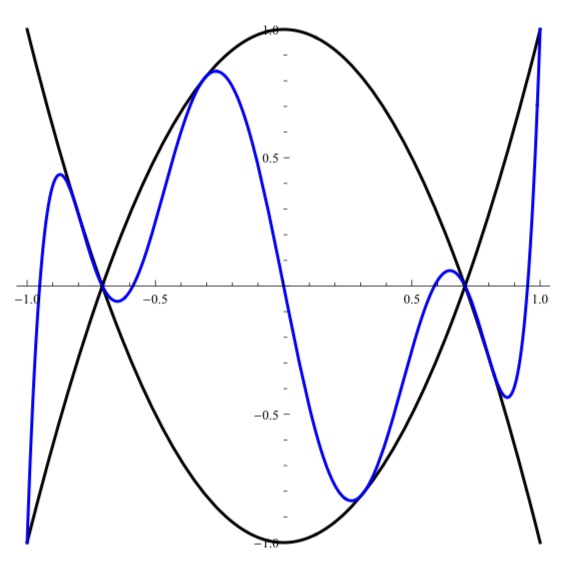}
         \caption{$\mu(x)=|1-2x^2|$}
         \label{fig:6}
     \end{subfigure}
        \caption{Examples of snake polynomials $\omega_\mu$.}
        \label{fig:snake}
\end{figure}
\begin{conjecture}[Geno Nikolov]
    If $\mu\geq 0$ is a continuous even convex function in $[-1,1]$, then the associated with $\mu$ snake polynomials have non-negative expansion in the Chebyshev polynomials of the first kind.
\end{conjecture}

To understand the motivation for these questions, let us recall some well-known inequalities for the derivatives of polynomials. Here $T_n(x)$ is the Chebyshev polynomial of the first kind of degree $n$, and $\|f\|=\max\limits_{x\in[-1,1]}|f(x)|$ is the usual uniform norm.
\begin{theorem}[V. A. Markov 1892]
    If $f\in \mathcal{P}_n$ and $\|f\|\leq 1$, then for $k=1,\dots,n$,  $$\|f^{(k)}\|\leq \|T^{(k)}_n\|.$$ The equality is attained only for $f=\pm T_n$.
\end{theorem}
\begin{theorem}[R. J. Duffin and A. C. Schaeffer 1941]
      If $f\in \mathcal{P}_n$ satisfies $|f\left(\cos(\nu\pi/n)\right)|\leq 1 $ for $\nu=0,\dots,n$, then for $k=1,\dots,n$,  $$\|f^{(k)}\|\leq \|T^{(k)}_n\|,$$ with equality only if $f=\pm T_n$.
\end{theorem}
The condition $\|f\|\leq 1$ means that $|f(x)|\leq 1$ for $x\in[-1,1].$ One can consider a general problem under restriction $|f(x)|\leq \mu(x)$, where $\mu$ is an arbitrary majorant. 

Given $n,k\in\mathbb{N}$, $1\leq k\leq n$, and a majorant $\mu\geq 0$, define
$$M_{k,\mu}=\sup\left\{\|p^{(k)}\|:p\in\mathcal{P}_n,|p(x)|\leq\mu(x),x\in[-1,1]\right\},$$
$$D^*_{k,\mu}=\sup\left\{\|p^{(k)}\|:p\in\mathcal{P}_n,|p(x)|\leq\mu(x),x\in\delta^*\right\}.$$

Clearly, $M_{k,\mu}\leq D^*_{k,\mu}$, and the results of V. A. Markov (1892) and of R. J. Duffin and A. C. Schaeffer (1941) can be read as: for $\mu\equiv 1$, we have  
$$M_{k,\mu}=D^*_{k,\mu}=\|T_n^{(k)}\|, \quad 1\leq k\leq n.$$
In other words, the snake polynomial $\omega_\mu=T_n$ for $\mu\equiv1$ is extremal for both Markov- and Duffin-Schaeffer- type inequalities. The next step is to study for which other majorants $\mu$ the snake polynomial $\omega_\mu$ satisfies 
$$M_{k,\mu}=D^*_{k,\mu}=\|\omega_\mu^{(k)}\|.$$ 
The motivation for the open problem comes from the following result.
\begin{theorem}[A. Shadrin and G. Nikolov, \cite{Nikolov2012,Nikolov2014}]
Given a majorant $\mu\in C[-1,1]$, $\mu\geq0$, let $Q_n$ be the associated with $\mu$ snake polynomial of degree $n$. If $Q_n$ admits non-negative or sign-alternating expansion in the Chebyshev polynomials of the first kind, then
$$M_{k,\mu}=D^*_{k,\mu},$$ 
    that is, $Q_n$ is extremal in the Markov- and Duffin-Schaeffer- type inequalities for polynomials with majorant $\mu$.
\end{theorem}

\section{\textbf{Chromatic numbers of spheres}
\\ Danila Cherkashin}
The \textit{chromatic number} $\chi( A)$  of a set  $A\subseteq\mathbb{R}^d$ is the minimal number of colours needed to colour all points of $A$ in such a way that any two points at the distance 1 have different colours. In the case $A=\mathbb{R}^2$, this is the well-known Hadwiger–Nelson problem. Our interest is in the case of the Euclidean
sphere $A=\mathbb{S}^{d-1}(R)\subseteq\mathbb{R}^d$ of radius $R$. It clear that $\chi(\mathbb{S}^{d-1}(R))=1$ for $0<R<1/2$ and $\chi(\mathbb{S}^{d-1}(1/2))=2$.

As shown in \cite{Raigorodskii2012OnTC}, for any fixed $R>1/2$,
the quantity $\chi(\mathbb{S}^{d-1}(R))$  grows exponentially. We consider the upper bound obtained by R. Prosanov in \cite{Prosanov2018}.

\begin{theorem}[R. Prosanov]
    For $R>\frac{\sqrt5}{2}$ we have
    $$\chi(\mathbb{S}^{d-1}(R))\leq \left(\sqrt{5-\frac{2}{R^2}+4\sqrt{1-\frac{5R^2-1}{4R^4}}}+o(1)\right)^d.$$
\end{theorem}

The proof consists of three steps:

\textbf{Step 1.} Fix $0<\phi(R)<\frac{\pi}{4}$. Consider a set $X\subset \mathbb{S}^{d-1}(R)$ of maximal cardinality s.t. for all $x_1,x_2\in X$ the spherical caps $C(x_1,\phi)$ and $C(x_2,\phi)$ do not intersect (i.e. $X$ is a kissing configuration). Consider the Voronoi tiling $\Psi$ of $\mathbb{S}^{d-1}(R)$ corresponding to the set $X$.

\textbf{Step 2.} Shrink every cell of $\Psi$ with ratio $\lambda(R)$. The union $\Psi^\prime$ of the smaller cells is colored in the first colour. The density of this set is $(\lambda(R)+o(1))^d.$

\textbf{Step 3.} Use greedy (or random) method to cover the sphere with differently colored translates of $\Psi^\prime$. We have $\chi(\mathbb{S}^{d-1}(R))\leq(\lambda(R)+o(1))^{-d}.$
\\
\\
Since Step 2 and Step 3 can be applied to (almost) any configuration from Step 1, the goal is to refine Step 1. 

\begin{question}[Danila Cherkashin]
Is it possible to find another configuration in Step 1 that would improve the upper bound for $\chi(\mathbb{S}^{d-1}(R))$?
\end{question}

\clearpage

\section{\textbf{Approximately Hadamard Matrices}\\Stefan Steinerberger}

Let us recall that a \textit{Hadamard matrix} is a rescaled orthogonal matrix $A\in\{-1,1\}^{n\times n}$ satisfying 
$$\|Ax\|_2=\sqrt{n}\|x\|_2,$$
for all $x\in\mathbb{R}^n$. It is known that if $A$ is a $n\times n$ Hadamard matrix and $n\geq 4$, then $n$ needs to be a multiple of $4$. This is a famous conjecture: whether this necessary condition is also sufficient. This leads us to a natural question about \textit{approximately Hadamard} matrices. Are there $A\in\{-1,1\}^{n\times n}$ satisfying 
$$c\sqrt{n}\|x\|_2 \leq \|Ax\|_2\leq C\sqrt{n}\|x\|_2,$$
so that $C/c$ is small? The answer is positive (Dong–Rudelson, \cite{Dong2023}). The construction is complex, and the proof itself requires nontrivial techniques. But there is an easy strengthening of the statement for \textit{circulant} matrices. 

\begin{theorem}[Steinerberger, \cite{Steinerberger2024}]
    There exist universal $0 < c < C < \infty$ such that for all $n \geq 1$, there exists a
circulant matrix $A \in\mathbb{R}^{n\times n}$ whose entries are $\pm 1$ such that
$$c\sqrt{n}\|x\|_2 \leq \|Ax\|_2\leq C\sqrt{n}\|x\|_2.$$
\end{theorem}
The result is constructive, using flat Littlewood polynomials, which is itself a difficult task. Thus, the goal is to find simpler ones. 
\begin{problem}[Stefan Steinerberger]
There should be many nice explicit constructions of approximately
Hadamard matrices with $C/c$ guaranteed to be small.
\end{problem}

\end{document}